\documentclass{article}
\usepackage[a4paper, left=2cm, right=2cm, top=2.5cm, bottom=1.5cm]{geometry}
\usepackage{subfiles}
\usepackage{mathtools}
\usepackage{amssymb}
\usepackage{graphicx,color}
\usepackage{amsthm}
\usepackage{amsmath}
\usepackage{enumerate}
\usepackage{tikz}
\usepackage{float}
\usepackage{hyperref}
\usepackage{rotating}
\usepackage{algorithm}
\usepackage{algpseudocode}
\usepackage{subcaption}
\newtheorem{theorem}{Theorem}
\newtheorem{prop}{Proposition}
\newtheorem{corollary}{Corollary}

\newtheorem{remark}{Remark}
\newtheorem{lemma}{Lemma}

\newtheorem{assumption}{Assumption}
\newtheorem{example}{Example}

\def\ba{\begin{array}}
  \def\ea{\end{array}}
\def\bi{\begin{itemize}}
  \def\ei{\end{itemize}}
\def\mR{\mathbb{R}}

\def\mE{\mathbb{E}}
\def\m1{1}

\let\oldhat\hat

\renewcommand{\hat}[1]{\oldhat{\mathbf{#1}}}

\providecommand{\com[2]}{\begin{tt}[#1: #2]\end{tt}}

\begin{document}

\title{Low complexity convergence rate bounds for push-sum algorithms with homogeneous correlation structure}

\author{Bal\'azs Gerencs\'er%
  \thanks{B. Gerencs\'er is with the HUN-REN Alfr\'ed R\'enyi Institute of Mathematics, Budapest, Hungary and the E\"otv\"os Lor\'and University, Department of Probability and Statistics, Budapest, Hungary, {\tt\small gerencser.balazs@renyi.hu} } %
  \and %
  Mikl\'os Kornyik%
  \thanks{M. Kornyik is with the HUN-REN Alfr\'ed R\'enyi Institute of Mathematics, Budapest, Hungary, {\tt\small kornyik.miklos@renyi.hu}} %
 }

\maketitle

\begin{abstract}
  The objective of this work is to establish an upper bound for the almost sure convergence rate for 
  a class of push-sum algorithms.
The current work extends the methods and results of the authors on a similar low-complexity bound on push-sum algorithms with some particular synchronous message passing schemes
and complements the general approach of Gerencsér and Gerencsér from 2022 providing an exact, but often less accessible description.

Furthermore, a parametric analysis is presented on the ``weight'' of the messages, which is found to be convex with an explicit expression for the gradient. This allows the fine-tuning of the algorithm used for improved efficiency.

Numerical results confirm the speedup in evaluating the computable bounds without deteriorating their performance, for a graph on 120 vertices the runtime drops by more than 4 orders of magnitude.
\end{abstract}

\section{Introduction}
\label{sec:intro}

The importance of the fundamental building block studied has been well understood for some while, the goal of distributed averaging of certain initial values along a network, while only relying on local communication and low computational overhead \cite{blondel2005convergence}, \cite{tsitsiklis:phd1984}. The main lines of research concern application of methods in various domains, others on extending the methodology and understanding of the processes. The current paper follows the latter goal.

In the baseline setting, usually some structure is assumed on the communication, in terms of the matrix describing the linear update of the vector of values. Most often it is required to be doubly stochastic, or even symmetric. This setting is well understood \cite{tahbaz2009consensus}, see the survey \cite{MAL-051} also guiding to applications, discussion and further references.

With applicability in mind, extensions to handle various network deficiencies have been studied. In particular, the interest for algorithms capable of handling asynchronous, directed communication appeared, which in turn caused the representing update matrix not to be doubly stochastic anymore, still with the unchanged goal to compute the exact average. To overcome this, the successful scheme of \emph{push-sum} was proposed \cite{kempe2003gossip}, also referred to as \emph{ratio consensus} \cite{hadjicostis2014average} with further variants developed such as \emph{weighted gossip} \cite{benezit2010weighted}.
The aim of these algorithms is the same but without imposing communication events to happen synchronously or consistently over the network.

	Other realistic communication challenges have been worthwhile research targets, the push-sum concept has been extended along various tracks to deal with such aspects, e.g., including packet loss \cite{hadjicostis2015robust} \cite{olshevsky2018fully} or delay \cite{hadjicostis2014average}
        or malicious agents intruding \cite{hadjicostis2022trustworthy}.
	Besides crafting variants, effort has been made to better understand the effect of such communication imperfections to the baseline push-sum algorithms. The error of the final consensus value in the presence of packet loss has been analyzed \cite{gerencser2018push} as it has been performed for classic (linear) gossip \cite{frasca2013large}, \cite{olfati2004consensus}.
        Thanks to the simple structure and objective of the algorithm, it also serves as the communication, synchronization layer for more complex efforts, e.g., for the spectral analysis of the actual network \cite{kempe2008decentralized} or for designing distributed optimization algorithms \cite{nedic2014distributed}, \cite{nedic2017achieving}.

        Recently there has been an increased interest to also handle situations when not only the communication connectivity is allowed to change, but also the participants in the network, leading to the study of open multi-agent networks. Consensus processes have been studied in this context, both for linear schemes \cite{monnoyer2022open} and also for push-sum protocols \cite{sawamura2024distributed}.    
	
	{\bf Related work.} An essential question to assess the usability and efficiency of such algorithms is understanding their asymptotic behavior, the convergence of the processes and the rate at which it takes place.
        In the original work \cite{kempe2003gossip}, exponential convergence has been shown for the reference push-sum scheme. At that stage, it has been a qualitative observation, without going after the exact convergence rate.
	
	For a more quantitative insight, an important step was proving a remarkable upper bound for the almost sure (a.s.) rate of convergence along an infinite subset of the timeline \cite{iutzeler2013analysis}. Another direction has been proposed \cite{rezaeinia2019push} to capture the finite time behavior, capable of stating explicit bounds, although far from being optimal.
    Via a much more general and abstract treatment, the asymptotic \emph{tight bound for the rate of a.s.\ convergence} has been identified \cite{gerencsr2019tight} through the gap of the Lyapunov exponents of the (random) update matrix series. Although this sheds light on the rate, the practical concern is that this Lyapunov gap is known to be \emph{incomputable} in several cases \cite{tsitsiklis1997lyapunov}. 
	
	Following up on the results above, it was possible to learn and combine from \cite{iutzeler2013analysis} and \cite{gerencsr2019tight}, eventually leading to an upper bound on the a.s.\ rate for the i.i.d.\ case \cite{gerencser2022computable}. The bounding expression is obtained by manipulating a single (random) update matrix, as such providing a \emph{computable quantity}. The bounds are numerically appealing, however for a graph on $N$ vertices, matrices of size $N^2\times N^2$ are encountered, thus \emph{increasing in complexity} quickly as $N$ grows.
        
        In view of this bottleneck, a bound with low computational complexity has been developed \cite{gkpushsumsynchr2024}, still for getting guarantees on the a.s.\ exponential convergence rate. This needed the restriction of focus to networks and communication schemes with high symmetries. The gain on numerical complexity is a speedup factor of at least $N$ for computing the bounds, but often more in practice.

	{\bf Contributions.} On one hand we extend the easily accessible convergence rate estimates of \cite{gkpushsumsynchr2024} for a broader class of messaging protocols. While still requiring a substantially symmetric setup, this is only in terms of the network and the correlation structure of the messaging, which is sufficient to deduce rate guarantees that are simple but perform well.
        Furthermore, we analyze the effect of the proportion of transmission, which is often set in some arbitrary way. It turns out that the rate bounds are convex in the proportion, and can be efficiently optimized (supported by structural understanding), leading to fine-tuning the push-sum algorithm in a quasi-optimal way. A range of numerical experiments illustrate the analytical results above.
        
	{\bf Layout.} The rest of the paper is structured as follows. In the next section we formally define the averaging processes considered, giving specific use case examples as well.
        In Section \ref{sec:polybounds} the first main result is presented, on easily accessible convergence rate guarantees.
        This is followed by the results on optimization considerations to adjust a given algorithm in Section \ref{sec:optimization}.
        Numerical demonstration of the quality of the bounds obtained and the effect of optimization is shown in Section \ref{sec:numerical}. A number of Appendices are included to spell out the technical details of our analysis.

\section{Model description}
\label{sec:model}
Let $(A(t))_{t\geq 1}$ be a sequence of i.i.d.\ (non-negative) column-stochastic matrices of size $N\times N$ and $x(0)\in \mR^N$ be the initial data. The initial auxiliary weight variable is initialized as $w(0)=\mathbf 1 \in \mR^N$. Push-sum algorithms, orginally introduced in \cite{kempe2003gossip}, are characterized by the following time-dependent equations:
\begin{align*}
	x(t) &= A(t)x(t-1),\\
	w(t) &=  A(t)w(t-1).
\end{align*}

The estimate of the average of the initial values at coordinate $i$ is given by $z_i(t) = x_i(t)/w_i(t)$.
Under some mild conditions this is known to converge almost surely to average consensus, i.e. $z_i(t)\to \bar x := N^{-1}\sum_{j}x_j(0)$ as $t\to \infty$ for all $i$, shown already for a range of cases in \cite{kempe2003gossip}.
Concerning the efficiency of this scheme, there is a tight asymptotic bound on the rate of almost sure convergence in \cite{gerencsr2019tight}. However, due to the computational difficulties of retrieving it, more accessible approximate solutions and upper bounds are quantities of high interest.

Let us write $A(t)$ as 
$$  A(t) = I + C(t) - D(t) $$
where $C(t)$ consists of the message passing terms with $c_{ij}\geq 0$, $\sum_i c_{ij} \leq 1$, and $c_{ii} = 0$ $\forall i,j$, while $D(t)$ denotes the corrective diagonal matrix so that $A(t)$ is column-stochastic. Note that these conditions imply $a_{ij}(t)\geq 0$ for all $i,j$. We remark that the term $D(t)-C(t)$ is also known as the Laplacian matrix from a graph theoretical perspective.

Let us prescribe the following relations concerning the elements of $C$.
\begin{align} \label{corr_rel}
	\mathbf 1^\top \mathbb EC&= q \sqrt u \mathbf 1^\top  \nonumber \\  
	\mathbb E[c_{ji}c_{lk}]&= \begin{cases}
		\alpha q^2  p_{ji}p_{lk} & i\neq k, \\
		\beta  q^2 p_{ji}p_{li} & i=k, j\neq l, \\
		(qr)^2 p_{ji} & i=k, j=l.
	\end{cases} 
\end{align}
Here $P=(p_{ij})$ is a column-stochastic matrix of non-negative elements and $u,q,r \in [0,1]$.
One can think of $P$ as the collection of normalized recipient distribution in each column and of $q$ as the average proportion that is sent per message.

For the current analysis, some assumptions have to be imposed on the transfer matrices $(A(t))_{t\geq 1}$, or equivalently on $(C(t))_{t \geq 1}$ as follows. 
\begin{assumption}
  \label{ass:C}
  The following are required for the communication events:
  \begin{enumerate}
    \item $(C(t))_{t\geq 1}$ forms an i.i.d.\ series of non-negative matrices;
      \item $C(1)$ satisfies the moment conditions given in \eqref{corr_rel};
    \item $C(1)$ is allowable, i.e.\ no row of $I+C(1)-D(1)$ can be fully 0;
    \item $C(1)$ does not have extremely small positive elements, more precisely for $m = \min\{C_{ij}(1) \mid C_{ij}(1)>0, i,j\in[N]\}$ we have $\mE\log^- m >-\infty$.
  \end{enumerate}
\end{assumption}

Furthermore, we need connectivity and a strong symmetry on the network where the averaging takes place:
\begin{assumption}
  \label{ass:symm}
  The following are required for the network descriptor $P$:
  \begin{enumerate}
  \item Connections are symmetric, i.e.\ $P=P^\top$;
\item The graph formed by the positive elements of $P$ is connected and $q\sqrt{u}<1$ to provide loops;
\item The network is transitive, i.e.\ for any $i,j$ vertices there exists a permutation matrix $\Pi$ so that $\Pi^{-1} P \Pi = P$ and the permutation maps $i$ to $j$.
  \end{enumerate}
\end{assumption}

To understand the motivation of \ref{corr_rel} let consider the following. Let us write the $(j,i)$ element of $C$ as
\begin{align*}
	c_{ji} & = Q_i U_i \frac{\chi_{ji}}{M_i}
\end{align*}
where $Q_i$ denotes the proportion sent by vertex $i$, $U_i \sim Bern(u)$ denotes the indicator of $i$ sending a message at all, $\chi_{ji}$ is the indicator of an interaction $i\to j$, while $M_i = \sum_{j} \chi_{ji}$ is basically the number of recipients. Note that in case vertex $i$ send messages to $M_i=m$ recipients, then the expected message rate per recipient will be $\tilde q/m$. We will see in this case that $q=\tilde q \sqrt u$. We make the following assumptions:
\begin{align*}
		\mathbb E[U_i] & = u, \\
	\mathbb E[Q_i|U_i]&= \tilde qU_i, \\
	\mathbb E[c_{ji} | U_i] &=  \tilde q  p_{ji}U_i,\\
	\mathbb E\left[\frac{\chi_{ji}}{M_i} | U_i\right] &= p_{ji}U_i, \\
	\mathbb E[Q_i^2|U_i]&= (\tilde qr)^2 U_i,
\end{align*} 
where $\sum_j p_{ji} = 1$ for all $i$, i.e.\ the matrix $(p_{ji})$ is column-stochastic. In the current context let us further assume that for any $i\neq k$ the quantities $Q_i, Q_k$ and $(\chi_{ji})_j, (\chi_{lk})_l$ are independent given $U_i,U_k$. Note that the condition $M_i>0$ is equivalent with $U_i = 1$, hence $\sigma(U_i) \subset \sigma(M_i)$, i.e. if a random variable $X$ is $M_i$-measurable then it is also $U_i$-measurable.
Let us now calculate the following quantities.
\begin{align*}
	\mathbb E[c_{ji}] &= \mathbb E[ qU_i\mathbb E[ \chi_{ji}/M_i | U_i]] = E[\tilde qU_i p_{ji}] =  \tilde q  u p_{ji}.
\end{align*}
For $i\neq k$ we have
\begin{align*}
	\mathbb E[c_{ji}c_{lk}] &= \mathbb E[\tilde q^2U_iU_k\mathbb E[\chi_{ji}/M_i|U_i] \mathbb E[\chi_{lk}/M_k|U_k]] = \tilde q^2 p_{ji}p_{lk} \mathbb E[U_iU_k] .
\end{align*}
Denote by $\alpha u := \mathbb E[U_iU_k]$, and obtain the trivial upper bound
$$ \alpha u := \mathbb E[U_i U_k] \leq \mathbb E[U_i] = u \Rightarrow 0\leq \alpha \leq 1.  $$
In case $i=k$ and $j\neq l$ the calculation reads
\begin{align*}
	\mathbb E[c_{ji}c_{li}] &= \mathbb E[\mathbb E[Q_i^2U_i \chi_{ji}/M_i \chi_{li}/M_i|U_i]]= \mathbb E[Q_i^2 U_i \mathbb E[\chi_{ji}\chi_{li}/M_i^2|U_i]] =: \beta u \tilde q^2  p_{ji}p_{li} \Rightarrow \beta\geq 0,
\end{align*}
while for $i=k$ and $j=l$, we obtain
\begin{align*}
	\mathbb E[c_{ji}^2] &= \mathbb E[Q_i^2 U_i \chi_{ji}/M_i^2] = \mathbb E[Q_i/M_i\mathbb E[Q_iU_i\chi_{ji}/M_i| M_i]]\\
	&= \mathbb E[ Q_i/M_i \mathbb E[c_{ji}|M_i]] = qp_{ji}\mathbb E[Q_i/M_i U_i] =: (r\tilde q)^2 p_{ji}u. 
\end{align*}
After summarizing the above we arrive at
\begin{align*}
	\mathbb E[c_{ji}c_{lk}]= \begin{cases}
		\alpha \tilde q^2 u p_{ji}p_{lk} & i\neq k, \\
		\beta  \tilde q^2 up_{ji}p_{li} & i=k, j\neq l, \\
		(\tilde qr)^2  u p_{ji} & i=k, j=l,
	\end{cases}
\end{align*}
which coincides with (\ref{corr_rel}) by setting $q: = \tilde q \sqrt u$.
Note that in this case
$$ \mathbb E[c_{ji}] = q\sqrt u p_{ji}.$$\\
The next examples represent the diversity of the model class.
\begin{example}\label{async_gossip}
	Consider a $d$-regular transitive graph $G=(V,E)$ and define the following dynamics: in each step a uniformly random vertex sends a message to a uniformly random neighbor with some globally fixed messaging rate $\tilde q\in(0,1)$. In this case $\alpha=\beta=0$, $r=1$ and $q=\tilde q \sqrt u$, $u = 1/N$ and $p_{ij} = d^{-1}\mathbb I((j,i)\in E)$.  
\end{example}

\begin{example}\label{ex_1}
	Let $G$ be a $d$-regular transitive graph and let us consider the following messaging protocol: in each step a randomly chosen vertex sends a message to a single random neighbor with probability $\rho $, or it broadcasts, i.e.\ sends a message to all of its neighbors, with probability $1-\rho$, and the total weight of the messages is given by $\tilde q$. The proportion sent can be specified, but in the simplest setting the total weight of the message(s) is equal in each case. This gives us the following:
	\begin{align*}
		\mathbb E[c_{ji}] & = \frac{\tilde q}{ N} \big(\rho d  + 1-\rho  \big) p_{ji}, \\
		\mathbb E[c_{ji}c_{lk}]&= \begin{cases}
			0 & i\neq k, \\
			\tilde q^2/N\ (1-\rho) p_{ji}p_{kl} & i=k, j\neq l, \\
			\tilde q^2/N\ \big(\rho d + (1-\rho)/d\big)p_{ji} & i=k, j=l,
		\end{cases}
	\end{align*}
and so $u=1/N$, $q = \tilde q \sqrt u$, $\beta = (1-\rho)$, $r^2 = \rho d + (1-\rho)/d$, $u=1/N$ and $ p_{ji} = \chi_{ji} /d $, where $\chi_{ji}$ is the indicator of the activity of edge $(i,j)$. 
\end{example}

\begin{example}\label{example:num}
	Given a $d$-regular transitive graph $G = (V,E)$, let $\chi_{ij}, (j,i)\in E$ denote i.i.d. random variables with Bernoulli($\rho$) ($\rho\in (0,1)$) distribution and let $$ c_{ji} =\tilde q \frac{\chi_{ji}}{\max\left(1,\sum\limits_{k: (i,k)\in E} \chi_{ki}\right)}. $$
	where $\tilde q \in (0,1)$ is the global fixed message rate for each node $j$. 
\end{example} The previous example is used to showcase our algorithm in Section \ref{sec:numerical}.

\begin{remark}
Although our results will not deal with them, we note that a natural extension to the scheme above is to specify node- and vertex-dependent coefficients as follows:
\begin{align} \label{corrstruct-general}	
	\mathbb E[c_{ji}] & = q_i\sqrt{u_i}p_{ji} \quad i\neq j, \nonumber \\
	\mathbb E[c_{ji}c_{lk}] &=  \begin{cases}  \alpha_{ijkl}q_i q_k p_{ji}p_{lk} & i\neq k ,\\ 
		\beta_{ijk}  q_i^2 p_{ji}p_{li}  & i=k, j\neq l, \\
		 (q_ir_i)^2 p_{ji} & i=k, j=l. 
	\end{cases}  
\end{align}  
Similar conditions are required for $P=(p_{ji})$ and $q_i,u_i,r_i \in [0,1]$ for all $i$.
\end{remark}

Let us outline further examples for this more general class. Note that Example
\ref{ex:indep-sync} fits in the stricter scheme in case $q_i$ are constant.
\begin{example} \label{ex:indep-sync}
	Each vertex sends a message independently to one of its neighbors, i.e. $\mathbb E[c_{ji}c_{lk}] = \mathbb E[c_{ji}]\mathbb E[c_{lk}] = u_iq_ip_{ji}u_kq_kp_{lk}. $
\end{example}
\begin{example}
	Given is a graph $G=(V,E)$, suppose there is a probability distribution $\nu$ given on the subset of vertices $2^V$. In each step a $\nu$ distributed random subset of vertices is chosen to be senders. Each of these vertices sends a message to a randomly chosen neighbor according to the distribution $(p_{ij})_i$ for sender $j$. The messages sent are independent conditioned on their senders. The overall proportions $q_i$ may be inhomogeneous. 
\end{example}

\section{Low-complexity rate bounds}
\label{sec:polybounds}

We will now state our first main result concerning the upper bound on the almost sure convergence rate of family of models described by \eqref{corr_rel}.

\begin{theorem}\label{generalprotocol} 
  Let Assumptions \ref{ass:C} and \ref{ass:symm} be in force for a homogeneous dynamics on $N$ vertices.
  Assume furthermore that $\forall i: p_{ii}=0$ and either $\forall i,j: p_{ij}\in \{0,c\}$ for some global $c\in (0,1)$ or $\beta =0$. Let $1=\lambda_1\geq \lambda_2 \geq \ldots \geq \lambda_N$ denote the real eigenvalues of $P$. Then for the rate of almost sure convergence we have
        \begin{equation}
          \label{eq:decayrate}
          \limsup_t \frac1t \max_i \log \bigg|\frac{x_i(t)}{w_i(t)} - \bar x\bigg|\leq \frac12 \log \xi_1    =: \frac 12 \gamma
        \end{equation}
	where $\xi_1$ denotes the largest root of the polynomial
	$$f(\xi)  = \prod_{i>1} (\xi - \Delta_i) \bigg(1- \frac{q^2}N \sum_{j>1} \frac{b_j}{\xi -\Delta_j}\bigg).  $$
	The constants $\Delta_j, b_j$ satisfy
	\begin{align*}
			\Delta_j &= 1-2q\sqrt u(1-\lambda_j) + \alpha q^2(1-\lambda_j)^2, \\
			b_j 
			& = (1-\lambda_j)\bigg((\beta-\alpha)(1-\lambda_j) - \beta c + 2r^2\bigg).
	\end{align*}

      \end{theorem}
As we can see, to obtain the bound above, simply the spectrum of a matrix of size $N\times N$ and a single root of a degree $N$ polynomial is required. This has theoretical worst case complexity of $O(N^3)$ \cite{dai2008convergence}. The numerical precision and practical computational cost is further discussed in Section \ref{sec:numerical}.\\
This theorem is a generalization of Theorem 2 in \cite{gkpushsumsynchr2024} and uses very similar methods, thus we will only focus on the novel elements.
Denoting $H(t) = A(t)A(t-1)\ldots A(1)$, the primary bound towards \eqref{eq:decayrate} is obtained as
\begin{equation}
  \begin{aligned} 
    \max_i \left|\frac{x_i(t)}{w_i(t)}-\bar x\right| &\leq \frac c {\min_i w_i(t)} ||x(0)||_2 ||H(t)(I-J)||_2\\ &\hspace{-1mm}\leq  \frac c{\min_iw_i(t)} ||x(0)||_2 ||H(t)(I-J)||_F \label{consensus_copy},
  \end{aligned}
\end{equation}
with some constant $c>0$ thanks to the equivalence of $l^2$ and Frobenius norms.
The first term on $w_i(t)$ is subexponential according to Lemma 10 of \cite{gerencser2022computable}, relying on Assumptions \ref{ass:C},\ref{ass:symm}. Note that in the cited lemma primitivity is needed, i.e.\ $H(t)$ eventually becoming elementwise positive, which is now provided by the connectivity assumption in terms of $P$ and requiring $q\sqrt{u}<1$, which ensures a positive diagonal of $\mE A(1)$.

The crucial last term is then handled using the map $\Phi(X) = \mathbb E [A(1)^\top X A(1)]$ and its adjoint $\Phi^*(X) = \mathbb E[A(1)XA(1)^\top]$ as
\begin{equation}
  \label{eq:errorL2bound}
  \| x(t) - \overline x \mathbf 1\|_2^2\leq \|(\Phi^*)^t(I-J))\|^2_F \|x(0)\|_2^2 \leq \|\Phi^*\|^t_{\mathcal x(0)} \|x(0)\|^2_2,
\end{equation}
where $\mathcal X_0 = \{X \in \mathbb R^{N\times N} : X=X^\top , XJ=0\}$.
Further details on the properties of $\Phi$ and the connections above are spelled out in Appendix \ref{sec:app_phi}. 
To proceed, we present some simple structural observations on $\Phi^*$.

\begin{lemma}\label{poly}
	Let $X$ be an arbitrary $N\times N$ matrix with constant diagonal and assume that $P$ is symmetric. Consider the general model (\ref{corrstruct-general}) with $\alpha_{ijkl}= \alpha$, $q_i=q$, $r_i=r$, and $\beta_{ijk} = 0$, or, in case $\beta_{ijk} = \beta \neq 0$, assume that $P$ consists of the elements $0,c$. Then $\Phi^*(X)$ is a degree 2 matrix polynomial of the variables $X,P,\omega(X)I$, where $\omega(X)$ denotes the common diagonal element of $X$.  
\end{lemma}
\noindent The proof can be found in the Appendix B.
\begin{corollary}\label{PX-commute}
  Suppose \eqref{corr_rel} is in force,
  and $P$ is symmetric and transitive.
  Then $X$ having constant diagonal $\omega$ and $XP=PX$ imply $P\Phi^*(X)=\Phi^*(X)P$.
\end{corollary}
For symmetric matrices the condition $XP=PX$, i.e.\ $X$ and $P$ commuting, implies that $P$ and $X$ have shared eigenvectors and there is a natural correspondence between their eigenvalues. We say that eigenvalues $\mu$ of $X$ and $\lambda$ of $P$ are in correspondence if they belong to the same shared eigenvector, i.e. $Pv = \lambda v$ and $Xv = \mu v$. This correspondence is well defined due to the orthogonality and hence linear independence of the eigenvectors.
\noindent The next proposition provides the backbone of our argument.
\begin{prop}\label{prop:recursion-eigenvalues}
  Let $X_0 = I-J$ and $X_{t+1} = \Phi^*(X_{t})$ and let Assumptions
  \ref{ass:C},\ref{ass:symm} be in force.
  Let us order the eigenvalues $\lambda_j$ (and the corresponding eigenvectors $v_j$) of $P$ in decreasing order, furthermore let $\mu_{t,j}$ denote the eigenvalue of $X_t$ corresponding to $v_j$ of $P$. Assume that either $\forall i,j:p_{ij}\in \{0,c\} $ for some global $c\in(0,1)$ or $\beta=0$. Then we have
	$$ \mu_{t+1,j} = \Delta_j \cdot \mu_{t,j} +q^2 \frac{b_j}{N} \sum_i \mu_{t,i} $$
	where $\Delta_j, b_j$ are defined in Theorem \ref{generalprotocol}.

\end{prop}
\noindent The proof is a series of quite straightforward but rather cumbersome calculations, its outline can be found in Appendix \ref{sec:app_phispelledout}.
 With the statements above in hand we have all the tools needed to prove Theorem \ref{generalprotocol}.
\begin{proof}[Proof of Theorem \ref{generalprotocol}]
	We can write
	$$  \mathrm{Tr}\ X_t = \sum_j \mu_{t,j},  $$
	and due to the Proposition \ref{prop:recursion-eigenvalues},
	$$ \mu_{t,j} = \Delta_j \cdot \mu_{t-1,j} + q^2\frac{b_j}N \sum_k \mu_{t-1,k}. $$
	Introducing the notations $\boldsymbol{\mu}_t = (\mu_{t,1},\ldots,\mu_{t,N})^\top$, ${\bf b}=(b_1,\ldots, b_N)^\top$ and $\boldsymbol{\Delta} = \sum_j\Delta_je_je_j^\top $  yields the compact form
	$$ \boldsymbol{\mu}_t = (\boldsymbol{\Delta} + q^2/N\ {\bf b}\mathbf 1^\top)\boldsymbol{\mu}_{t-1}. $$
	This is the same setup as in our previous work \cite[Proposition 3]{gkpushsumsynchr2024} but with different $\boldsymbol\Delta$ and $ {\bf  b}$ coefficients, in which we have shown that 
	\[\lim_t  \frac1t \max_j \log \mu_{t,j} \leq (\log \xi_1)/2 ,\]
	where $\xi_1$ denotes the largest root of the characteristic polynomial
	\begin{equation}
		\prod_{i>1} (\xi - \Delta_i) \bigg(1- \frac{q^2}N \sum_{j>1} \frac{b_j}{\xi -\Delta_j}\bigg). \label{RatePolynomial}
	\end{equation}  
        Given that the claim does not depend on the actual $\boldsymbol\Delta, {\bf  b}$ used we conclude.
\end{proof}

\section{Optimizing transmission intensity for rates}
\label{sec:optimization}

In what follows we analyze the dependence on $q$ of the bound obtained in Theorem \ref{generalprotocol}, which will be helpful for algorithmic tuning. Even though optimizing an upper bound on the rate does not apriori provide the optimal rate, we will see this is an easily accessible quantity with good performance, see also Section \ref{sec:numerical} for numerical demonstration.

Establishing the convexity of $q\to\xi_1(q)$ is interesting on it own,
and opens a wide range of possibilities for optimization.

\begin{theorem}\label{ConvThm}
  Assume the same conditions as in Theorem \ref{generalprotocol}. Then
the function $[0,1]\ni q \mapsto \xi_1(q)$ is convex.
\end{theorem}
\begin{proof} 
	By definition $\xi_1(q)$ is the largest solution to $f(\xi, q) = 0 $ in $\xi$. The latter can be written as
	\begin{equation} \label{ConvThmEq1}
\prod_{j>1} (\xi - \Delta_j) \bigg( 1- \frac{q^2}{N} \sum_{i>1} \frac{b_i}{\xi- \Delta_i}\bigg) = 0. 
	\end{equation} 
	According to Proposition \ref{LargestRootIneq} detailed in Appendix \ref{app:optimal} we know that $\xi_1(q) > \Delta_{j_0}$, the largest entry of $\mathbf \Delta$, meaning that the first term on the left side of \eqref{ConvThmEq1} is non-zero, therefore 
	$$ \frac{q^2}{N}\sum_{j>1} \frac{b_j}{\xi - \Delta_j} = 1 $$
	must hold for $\xi=\xi_1(q)$. Due to Corollary \ref{PosDefHessian} (see Appendix \ref{app:optimal}) we know that each function $h_j:(\xi,q) \mapsto q^2 b_j / (\xi-\Delta_j)$ is convex if $\alpha\geq 0$ and $\xi \geq \Delta_{j}$, meaning that the sum of these functions is convex if $\alpha \geq0$ and $\xi \geq \Delta_{j_0}$. This implies that the level sets of $\sum_j h_j$ are also convex, hence $\xi_1(q)$, as a solution to $\sum_j h_j = N$ as well.
      \end{proof}
In addition to convexity, the gradient can be expressed in close form once $\xi_1$ is evaluated, which is a useful additional quantity for the application of a range of optimizers:

      \begin{prop}
  \label{prop:lines_xidiff}
	Assume that the conditions of Theorem \ref{ConvThm} hold, furthermore $\lambda_j>-1$ for all $j$. Then the rate $\xi_1(q)$ satisfies $$1+\xi_1'(0)q \leq \xi_1(q) \leq 1+\xi_1'(1)(q-1) \qquad \forall q\in(0,1)$$
	where
	\begin{align*}
		\xi_1'(0) &= -2(1-\lambda_2), \\
		\xi_1'(1) &=2\  \frac{\sum_{j>1}(1-\lambda_j)b_j/(2(1-\lambda_j) - \alpha(1-\lambda_j)^2)^2 }{\sum_{j>1} b_j/(2(1-\lambda_j) - \alpha(1-\lambda_j)^2)^2 },
	\end{align*}
	and generally
	\begin{equation}
          \xi_1'(q) = \frac{2q \sum_{j>1} b_j/(\xi_1(q)-\Delta_j) + q^2 \sum_{j>1} b_j \partial_q \Delta_j / (\xi_1(q)-\Delta_j)^2}{q^2\sum_{j>1} b_j/(\xi_1(q)-\Delta_j)^2}. \label{eq:xidiff}
          \end{equation}
\end{prop}
\noindent
For the proof of Proposition \ref{prop:lines_xidiff} see Appendix \ref{app:optimal}.

Once a convex function is available together with its gradient, standard convex optimization techniques can be applied for minimization \cite{nesterov2018lectures}. Numerical experiments, case studies are presented in the next section, analyzing the computation cost of optimization and the efficiency gain obtained.

\section{Numerical experiments}
\label{sec:numerical}

As we mentioned in Section \ref{sec:model} we will showcase our algorithm on Example \ref{example:num}. In order to realize the dynamics let $\chi_{ij}, (j,i)\in E$ be independent variables of $Bern(\rho)$, indicators of the activity of edge $(j,i)$ at a given time instance, $\rho,\tilde q \in (0,1)$, and $G = (V,E)$ a $d$-regular, transitive graph. Set
$$c_{ij} = \tilde q\frac{\chi_{ji}}{\max\left(1,\sum\limits_{k: (i,k)\in E} \chi_{ki}\right)}.$$  Straightforward calculation shows
\begin{align*}
	\mathbb E[c_{ji}] & = \tilde q (1-(1-\rho)^d) p_{ji} \\
	\mE [c_{ji}c_{lk}]& = \begin{dcases}
		\tilde q^2 (1-(1-\rho)^d)^2 p_{ji}p_{lk}, & \mbox{ for } i\neq k, \\
		\tilde q^2d^2 \sum_{m=2}^d \frac{1}{m^2} \binom{d-2}{m-2} \rho^m (1-\rho)^{d-m} p_{ji}p_{lk}  & \mbox{ for } i=k, j\neq l, \\
		\tilde q^2 d\sum_{m=1}^d \frac{1}{m^2} \binom{d-1}{m-1} \rho^m (1-\rho)^{d-m} p_{ji} & \mbox{otherwise.}
	\end{dcases}
\end{align*}
where $p_{ij} = 1/d$ if $(j,i) \in E$ else 0. 
By setting $u = (1-(1-\rho)^d)$, and $q = \tilde q \sqrt u $ the previous relations assume the same form described in \eqref{corrstruct-general} with $\alpha = u$, $\beta = \frac{d^2}{u} \sum_{m=2}^d \frac{1}{m^2}\binom{d-2}{m-2} \rho^m(1-\rho)^{d-m} $ and $r^2 =  \frac{d}{u}\sum_{m=1}^d \frac1{m^2} \binom{d-1}{m-1}\cdot $ $ \rho^m(1-\rho)^{d-m} .$

We will now illustrate both the precision and the computation complexity of our method given in Theorem \ref{generalprotocol}. For the choice of graphs, on one hand we use random Cayley graphs of the symmetric group $S_k$, i.e.\ taking random permutations on $k=5$ elements, then use the generated graph on the $5!=120$ vertices. We also identify an instance when only the alternating group $A_5$ is generated, providing a transitive graph on $60$ vertices.
On the other hand, to explore in a less expanding context, we use a grid of dimensions $5\times 6\times 7$, on the torus for transitivity, providing 210 vertices.

For comparing computational efficiency of our method, we present in Table \ref{tab:benchmark} below the computation time in seconds $T(\cdot)$ of our bound $\gamma/2$ and of $\eta_2/2:= \log \rho\big(\mE (A(1)^{\otimes 2}) (I-J)^{\otimes 2}\big)/2$ given in \cite{gerencser2022computable}, which is applicable in more general settings, but relies on higher dimensional objects.
As expected, the algorithm developed is considerably faster, by a factor approximately $N^2$.
Computations have been performed on AMD EPYC 7643 CPUs using the Julia (v1.10) computing language \cite{Julia2017}.

\begin{table}[h] 
  \small
	\centering
	\begin{tabular}{|c|c|c|c|} \hline
		$N$ & $60$ & $120$ & $210$  \\ \hline 
		$T(\gamma)$	& 0.00015 & 0.00033 & 0.0011\\ \hline 
		$T(\eta_2)$ & 0.18516 & 9.12899	& 136.40832 \\ \hline 
		$\log_N\big(T(\eta_2)/T(\gamma)\big)$ & 1.8143 & 2.136  & 2.193 \\ \hline 
	\end{tabular}
	\caption{Runtime in seconds comparing the computation of $\gamma$ and $\eta_2$.}
        \label{tab:benchmark}
      \end{table}
For precision, we compare with an approximation of the true convergence rate via simulating the dynamics, evaluating
\[
  \zeta := \frac{1}{t}\log \bigg\|\frac{1}{\sqrt{M}} diag(w(t))^{-1} H(t)   X^M(0) \bigg\|_F,
\]
where $X^M(0)$ is an $ N\times M$ matrix of uniform random independent columns in $\mathbf 1^\perp$ with unit norm. This is an extension from using a single initial vector $x(0)$ for evaluation, rather jointly a collection of them, reaching the principal rate with high probability. We chose $M=\lfloor \sqrt{N} \rfloor$ and $t=20000$.

\begin{figure}[h!]
	
	\centering 
	\begin{subfigure}{0.4\textwidth}
		\centering 
		
		\includegraphics[width=\textwidth]{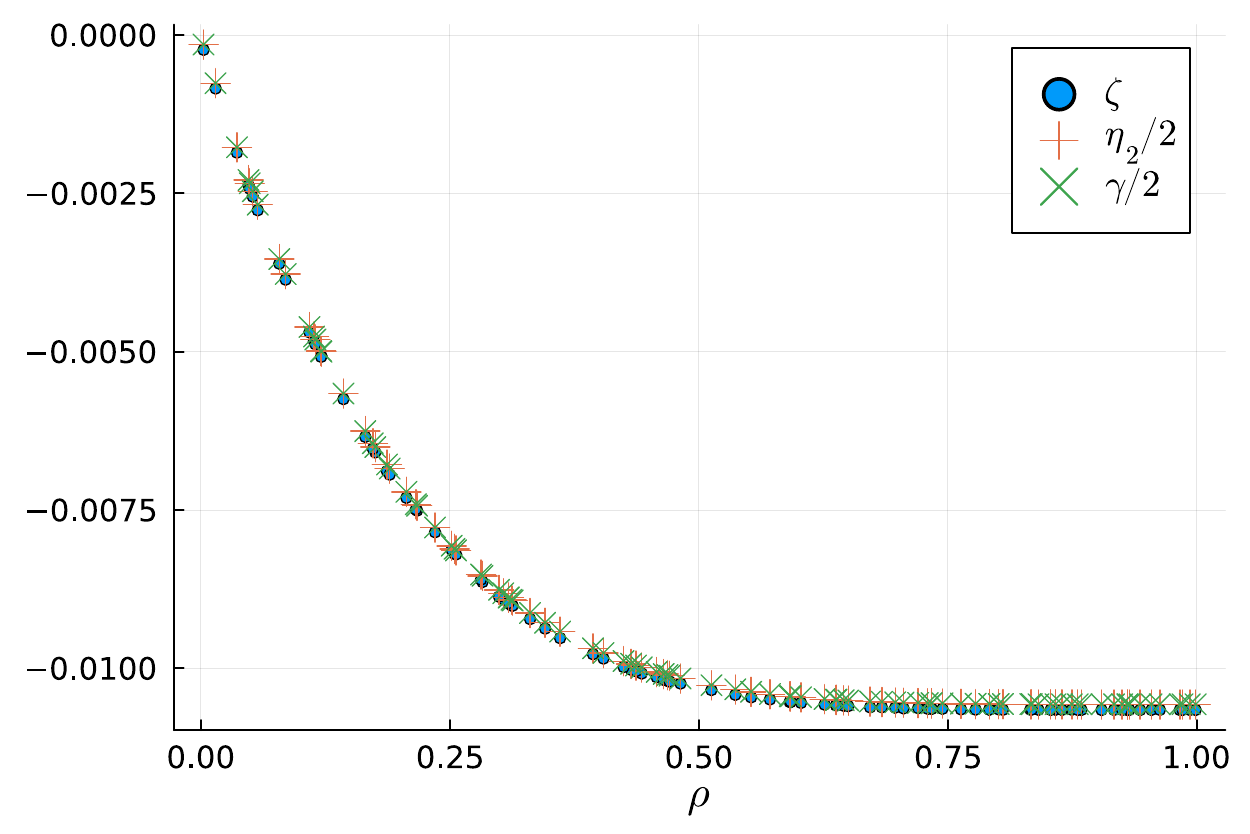}
		\caption{Fixed message proportion $\tilde q = 0.1$ and 100 random values of $\rho$.}
		\label{fig::numerics-a}
	\end{subfigure} \hfill
	\begin{subfigure}{0.4\textwidth}
		\centering 	
		\includegraphics[width=\textwidth]{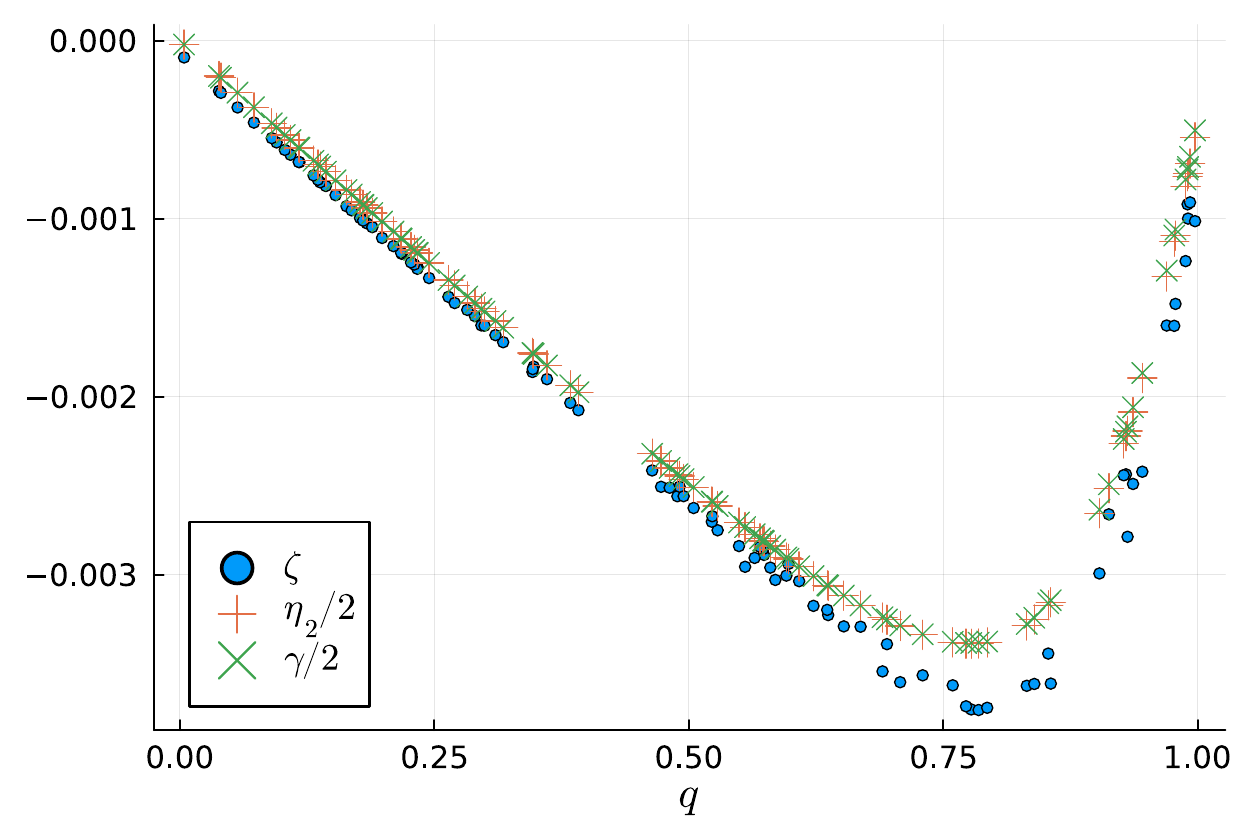}
		\caption{Fixed message probability $\rho=0.01$ and 100 random samples of the parameter $\tilde q$.}
			\label{fig::numerics-b}
	\end{subfigure}
\caption{Approximate rate by simulation ($\zeta$) along with bound using tensors ($\eta_2/2$) and of Theorem \ref{generalprotocol} ($\gamma/2$). }
\label{fig::numerics}
\end{figure}

Together with the error of the bounds we show also the interdependence of the parameters, 
Figures \ref{fig::numerics-a} and \ref{fig::numerics-b} show $\rho$ and $\tilde q$ dependence, respectively, of the different bounds $\gamma/2,\eta_2$ and the (approximate) Lyapunov exponent $\zeta$. According to these figures, the bounds provide reliable estimates, with the gap from the true rate being sometimes an order of magnitude smaller than the rate itself, but mostly much less.

\begin{figure}[h] 
  \centering 
  \includegraphics[width=0.6\textwidth]{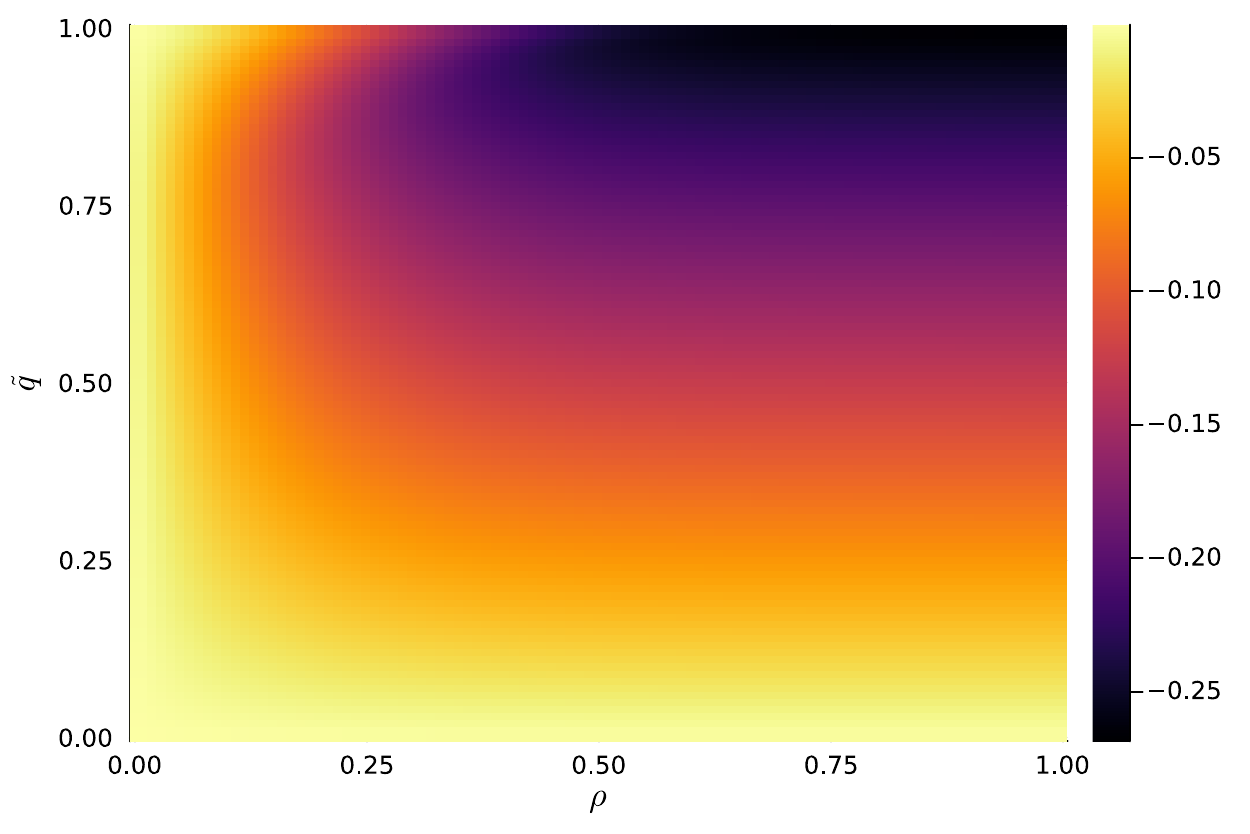}
  \caption{Heatmap showing our bound $\gamma/2$ on a transitive graph of 720 nodes, depending on the parameters $\tilde q$ and $\rho$.}
  \label{fig::heatmap}	
  
\end{figure}

In Figure \ref{fig::heatmap} we show the joint dependence on $\rho,\tilde q$ of our bound $\gamma/2$ via a heatmap. Again a random Cayley graph is used on $S_6$, a graph on 720 nodes and degree.
An interesting phase transition can be observed when considering tuning $\tilde q$ for fixed $\rho$: for small values there is a non-trivial optimal $\tilde q$ while for larger $\rho$ it is better to set $\tilde q = 1$. 
In general it seems preferable to increase the parameters to $\rho=\tilde q=1$, however if a cost of messaging is included in some way, the optimal choice might be different.

\section*{Acknowledgements}
The research was supported by NRDI (National Research, Development and Innovation Office) grant KKP 137490.\\
Miklós Kornyik was partially supported by MILAB. \\
No AI was used in preparing this work.

\bibliographystyle{siam}
\bibliography{generic,pushsum}

\begin{thebibliography}{10}

\bibitem{benezit2010weighted}
{\sc F.~B{\'e}n{\'e}zit, V.~Blondel, P.~Thiran, J.~N. Tsitsiklis, and
  M.~Vetterli}, {\em Weighted gossip: Distributed averaging using non-doubly
  stochastic matrices}, in {P}roceedings of 2010 {IEEE} {I}nternational
  {S}ymposium on {I}nformation {T}heory ({ISIT}), 2010, pp.~1753--1757.

\bibitem{Julia2017}
{\sc J.~Bezanson, A.~Edelman, S.~Karpinski, and V.~B. Shah}, {\em Julia: A
  fresh approach to numerical computing}, SIAM {R}eview, 59 (2017), pp.~65--98.

\bibitem{blondel2005convergence}
{\sc V.~Blondel, J.~M. Hendrickx, A.~Olshevsky, and J.~N. Tsitsiklis}, {\em
  Convergence in multiagent coordination, consensus, and flocking}, in 44th
  IEEE Conference on Decision and Control, 2005, pp.~2996--3000.

\bibitem{dai2008convergence}
{\sc X.~Dai, J.~Xu, and A.~Zhou}, {\em Convergence and optimal complexity of
  adaptive finite element eigenvalue computations}, Numerische Mathematik, 110
  (2008), pp.~313--355.

\bibitem{frasca2013large}
{\sc P.~Frasca and J.~M. Hendrickx}, {\em Large network consensus is robust to
  packet losses and interferences}, in 2013 European Control Conference (ECC),
  IEEE, 2013, pp.~1782--1787.

\bibitem{gerencser2022computable}
{\sc B.~Gerencs{\'e}r}, {\em Computable convergence rate bound for ratio
  consensus algorithms}, IEEE Control Systems Letters, 6 (2022),
  pp.~3307--3312.

\bibitem{gerencsr2019tight}
{\sc B.~Gerencs{\'e}r and L.~Gerencs{\'e}r}, {\em Tight bounds on the
  convergence rate of generalized ratio consensus algorithms}, IEEE
  Transactions on Automatic Control, 67 (2022), pp.~1669--1684.

\bibitem{gerencser2018push}
{\sc B.~Gerencs{\'e}r and J.~M. Hendrickx}, {\em Push sum with transmission
  failures}, IEEE Transactions on Automatic Control, 64 (2018), pp.~1019--1033.

\bibitem{gkpushsumsynchr2024}
{\sc B.~Gerencsér and M.~Kornyik}, {\em Low complexity convergence rate bounds
  for the synchronous gossip subclass of push-sum algorithms}, IEEE Control
  Systems Letters, 8 (2024), pp.~1283--1288.

\bibitem{hadjicostis2014average}
{\sc C.~N. Hadjicostis and T.~Charalambous}, {\em Average consensus in the
  presence of delays in directed graph topologies}, IEEE Transactions on
  Automatic Control, 59 (2014), pp.~763--768.

\bibitem{hadjicostis2022trustworthy}
{\sc C.~N. Hadjicostis and A.~D. Dom{\'\i}nguez-Garc{\'\i}a}, {\em Trustworthy
  distributed average consensus}, in 2022 IEEE 61st Conference on Decision and
  Control (CDC), IEEE, 2022, pp.~7403--7408.

\bibitem{hadjicostis2015robust}
{\sc C.~N. Hadjicostis, N.~H. Vaidya, and A.~D. Dom{\'\i}nguez-Garc{\'\i}a},
  {\em Robust distributed average consensus via exchange of running sums}, IEEE
  Transactions on Automatic Control, 61 (2015), pp.~1492--1507.

\bibitem{iutzeler2013analysis}
{\sc F.~Iutzeler, P.~Ciblat, and W.~Hachem}, {\em Analysis of sum-weight-like
  algorithms for averaging in wireless sensor networks}, IEEE Transactions on
  Signal Processing, 61 (2013), pp.~2802--2814.

\bibitem{kempe2003gossip}
{\sc D.~Kempe, A.~Dobra, and J.~Gehrke}, {\em Gossip-based computation of
  aggregate information}, in Proceedings of 44th Annual IEEE Symposium on
  Foundations of Computer Science, 2003, pp.~482--491.

\bibitem{kempe2008decentralized}
{\sc D.~Kempe and F.~McSherry}, {\em A decentralized algorithm for spectral
  analysis}, Journal of Computer and System Sciences, 74 (2008), pp.~70--83.

\bibitem{monnoyer2022open}
{\sc C.~Monnoyer de Galland~de Carni{\`e}res}, {\em Open multi-agent systems:
  representation, limitations and decentralized optimization}, PhD thesis,
  UCL-Universit{\'e} Catholique de Louvain, 2022.

\bibitem{nedic2014distributed}
{\sc A.~Nedi{\'c} and A.~Olshevsky}, {\em Distributed optimization over
  time-varying directed graphs}, IEEE Transactions on Automatic Control, 60
  (2014), pp.~601--615.

\bibitem{nedic2017achieving}
{\sc A.~Nedi{\'c}, A.~Olshevsky, and W.~Shi}, {\em Achieving geometric
  convergence for distributed optimization over time-varying graphs}, SIAM
  Journal on Optimization, 27 (2017), pp.~2597--2633.

\bibitem{nesterov2018lectures}
{\sc Y.~Nesterov}, {\em Lectures on convex optimization}, vol.~137, Springer,
  2018.

\bibitem{olfati2004consensus}
{\sc R.~Olfati-Saber and R.~M. Murray}, {\em Consensus problems in networks of
  agents with switching topology and time-delays}, IEEE Transactions on
  Automatic Control, 49 (2004), pp.~1520--1533.

\bibitem{olshevsky2018fully}
{\sc A.~Olshevsky, I.~C. Paschalidis, and A.~Spiridonoff}, {\em Fully
  asynchronous push-sum with growing intercommunication intervals}, in 2018
  Annual American Control Conference (ACC), IEEE, 2018, pp.~591--596.

\bibitem{rezaeinia2019push}
{\sc P.~Rezaeinia, B.~Gharesifard, T.~Linder, and B.~Touri}, {\em Push-sum on
  random graphs: Almost sure convergence and convergence rate}, IEEE
  Transactions on Automatic Control, 65 (2019), pp.~1295--1302.

\bibitem{sawamura2024distributed}
{\sc R.~Sawamura, N.~Hayashi, and M.~Inuiguchi}, {\em A distributed primal-dual
  push-sum algorithm on open multiagent networks}, IEEE Transactions on
  Automatic Control,  (2024).

\bibitem{MAL-051}
{\sc A.~H. Sayed}, {\em Adaptation, learning, and optimization over networks},
  Foundations and Trends® in Machine Learning, 7 (2014), pp.~311--801.

\bibitem{tahbaz2009consensus}
{\sc A.~Tahbaz-Salehi and A.~Jadbabaie}, {\em Consensus over ergodic stationary
  graph processes}, IEEE Transactions on Automatic Control, 55 (2009),
  pp.~225--230.

\bibitem{tsitsiklis:phd1984}
{\sc J.~N. Tsitsiklis}, {\em Problems in decentralized decision making and
  computation}, PhD thesis, Massachusetts Institute of Technology, 1984.

\bibitem{tsitsiklis1997lyapunov}
{\sc J.~N. Tsitsiklis and V.~Blondel}, {\em The {L}yapunov exponent and joint
  spectral radius of pairs of matrices are hard—when not impossible—to
  compute and to approximate}, Mathematics of Control, Signals and Systems, 10
  (1997), pp.~31--40.

\end{thebibliography}

\appendix
\section*{Appendix}

\section{Technical bounding details}
\label{sec:app_phi}

The following reasoning along with the provided calculations were presented in \cite{gkpushsumsynchr2024}, but for the sake of convenience and better understanding we will include it here as well.

Note that the elements $w_i(t)$ stay all positive as $A(t)$ is allowable.
Using the notations $H(t)=A(t)A(t-1)\cdots A(1)$, $J=\mathbf 1 \mathbf 1^\top /N$ we can write
\begin{align*}  x(t) - \bar x w(t)  &=  H(t)x(0)- \bar x w(t) \\ 
	&= (H(t)J + H(t)(I-J))x(0) - \bar x w(t) \\ 
	&=H(t)(I-J)x(0).
\end{align*}
After dividing by $w_i(t)$ the error in question can be bounded, with some constant $c>0$ for the equivalence of norms, as
\begin{equation}
  \begin{aligned} 
    \max_i \left|\frac{x_i(t)}{w_i(t)}-\bar x\right| &\leq \frac c {\min_i w_i(t)} ||x(0)||_2 ||H(t)(I-J)||_2\\ &\hspace{-1mm}\leq  \frac c{\min_iw_i(t)} ||x(0)||_2 ||H(t)(I-J)||_F \label{consensus}.
  \end{aligned}
\end{equation}
Now to get an insight of the almost sure convergence rate, $1/\min_iw_i(t)$ can be seen to be subexponential via \cite[Lemma 10]{gerencser2022computable}, relying on Assumption \ref{ass:C} and the connectivity requirement of Assumption \ref{ass:symm}. Therefore the decay rate will be controlled by $||H(t)(I-J)||_F$. A recursive relation will be established, let us first look at the second moment:
\begin{align*}
  \mathbb E||H(t)(I-J)||_F^2 &= \mathbb E\ \mathrm{Tr}\{(I-J)H(t)^\top H(t) (I-J)\} \\ & = \mathrm{Tr}\{\mathbb E [H(t)(I-J)H(t)^\top ]\}\\
                             &= \sum_j \mu_{t,j} \leq N \max_j \mu_{t,j},
\end{align*}
with $\mu_{t,j}$ denoting the $j^{th}$ eigenvalue of $X_t:=\mathbb E [H(t)(I-J)H(t)^\top ]$.
Considering the expression above, even for an arbitrary matrix $Y$ independent of $A(t)$, we have 
\begin{align}
	\label{phi_patter}
	\mathbb E[H(t) Y H(t)^\top] & = \mathbb E \big[\mathbb E[A(t)  H(t-1) Y H(t-1)^\top A(t)^\top \mid H(t-1) ] \big] \nonumber \\& = \mathbb E[\mathbb E_{A(t)}\big.  [A(t) X A(t)^\top ]\big|_{X= H(t-1)YH(t-1)^\top}].
\end{align}
The relation can be expressed using the linear operator $\Phi: \mathbb R^{N\times N} \to \mathbb R^{N\times N}$ acting on matrices:
$$
\Phi(X):=\mathbb E[AXA^\top ].
$$

With this notation $X_t$ will satisfy the recursion $X_0 = I-J$ and $X_t = \Phi(X_{t-1})$. In order to get a handle on (\ref{consensus}) we will need to conduct a thorough analysis of $\Phi$.

 \section{Discussion of $\Phi$}
 \label{sec:app_phispelledout}
 
For the sake of simplicity will omit $t$ from $A(t)$, $C(t)$ and $D(t)$. Let us revise the definition of $\Phi^*$ and perform some basic calculations. 
\begin{align*}
	\Phi^*(X) &= \mathbb E [AXA^\top] =  X - \mathbb E[(D-C)X] - \mathbb E[X(D-C)^\top] + \mathbb E [(D-C)X(D-C)^\top]\\
	&=X-(I-P)QX - XQ(I-P^\top) + \mathbb E [(D-C)X(D-C)^\top].
\end{align*}
where we now permit the cases that the message rate and the second moment may vary from node to node, i.e. given are $q_i,r_i \in (0,1)$ for each $i$. Hence $Q,R$ denote the diagonal matrices with elements $q_i$ and $r_i$ resp. 
It is obvious that the only non-trivial term is $\mathbb E[(D-C)X(D-C)^\top]$, hence in order to prove Proposition \ref{prop:recursion-eigenvalues}

we will need to get a handle on this term.
For the formal treatment, let us define the operator of extracting the diagonal:
$$ \Psi : \mathbb R^{N\times N} \ni X \mapsto (x_{11}, x_{22}, \ldots, x_{NN}) \in \mathbb R^N  $$
and its pseudo-inverse
$$ \Psi^- : \mathbb R^N \ni x \mapsto \mathrm{diag}(x_{1},\ldots, x_{N}) \in \mathbb R ^{N\times N}.$$

\begin{lemma}
Suppose that the assumptions of Lemma \ref{poly} hold, furthermore let $X^0:= X - \Psi^- \Psi X $. Then
	\begin{align*}
		\mathbb E[DXC^\top] &= \alpha QX^0QP^\top + \beta Q\Psi^-\Psi(X) Q P^\top  - \beta Q\Psi^-\Psi(X) Q(P\odot P)^\top  + QR\Psi^-\Psi(X) RQ P^\top, \\
		\mathbb E[DXD]&=  \alpha QX^0Q + \beta \Psi^-\Psi(X)Q^2(I-\Psi^- (P^\top \odot P^\top \ \mathbf 1 )) + QR\Psi^-\Psi( X) RQ, \\
		\mathbb E[C X C^\top]&= \alpha PQX^0QP^\top + \beta PQ \Psi^-\Psi (X) QP^\top - \beta \Psi^-(P\odot P\ \Psi(Q^2X)) + \Psi^-(P\ \Psi((QR)^2X)).
	\end{align*}
\end{lemma}

\begin{proof}
	
	\begin{enumerate}
		\item First let us calculate
		\begin{align*}
			(\mathbb E[DXC^\top])_{ij} &= \sum_{k} \mathbb E[d_{ii}x_{ik}c_{jk}] = \sum_{k} x_{ik} \sum_{l:l\neq i} \mathbb E[c_{li} c_{jk}]  \\
			& = \sum_{k:k\neq i} x_{ik} \sum_{l:l\neq i} \alpha q_iq_k p_{li}p_{jk}   +  x_{ii} \sum_{l:l\neq i,j} \beta q_i^2 p_{li}p_{ji} + x_{ii} q_i^2r_i^2 p_{ji} \\
			& = \alpha \sum_{k} x_{ik}q_iq_kp_{jk} - \alpha x_{ii} q_i^2p_{ji} + \beta x_{ii} q_i^2 p_{ji} - \beta x_{ii} q_i^2 p_{ji}^2 + x_{ii}q_i^2r_i^2 p_{ji}
		\end{align*}
		whence
		$$ \mathbb E[DXC^\top] = \alpha QX^0QP^\top + \beta Q\Psi^-\Psi(X) Q P^\top  - \beta Q\Psi^-\Psi(X) Q(P\odot P)^\top  + QR\Psi^-\Psi(X) RQ P^\top. $$
		\item Next let us compute
		\begin{align*} \mathbb E[(DXD)_{ij}]& = \mathbb E[d_{ii} x_{ij} d_{jj}] = \sum_{k:k\neq i} \sum_{l:l\neq j} x_{ij}E[c_{ki}c_{lj}] \\
			& = \delta_{ij} x_{ii}\big( \sum_{\substack{k,l:k\neq l \\ k\neq i \neq l}} \beta q_i^2p_{ki}p_{li} + \sum_{k:k\neq i} q_i^2r_i^2 p_{ki} \big) + (1-\delta_{ij})x_{ij} \sum_{\substack{k:k\neq i\\l:l\neq j }}\alpha q_iq_j p_{ki}p_{li}\\
			& =\alpha x_{ij}q_i q_j + \delta_{ij}(-\alpha q_i^2 + \beta q_i^2 - \beta \sum_{k\neq i}  q_i^2 p_{ki}^2 + q_i^2r_i^2)x_{ii}
		\end{align*}
		and so
		$$ \mathbb E[DXD] = \alpha QX^0Q + \beta \Psi^-\Psi(X)Q^2(I-\Psi^- (P^\top \odot P^\top \ \mathbf 1 )) + QR\Psi^-\Psi X RQ  .  $$
		\item Last we have
		\begin{align*}
			(\mathbb E[CXC^\top])_{ij}&=\sum_{k,l} \mathbb E [c_{ik}x_{kl}c_{jl}] \\
			&= \sum_{k\neq l} \alpha x_{kl} q_k q_l p_{ik} p_{jl} + \sum_k \mathbb E [c_{ik}x_{kk}c_{jk}] \\
			& = \alpha \sum_{k,l} x_{kl} q_k q_l p_{ik} p_{jk} - \alpha \sum_k x_{kk} q_k^2 p_{ik}p_{jl} + \sum_k \mathbb E [c_{ik} x_{kk} c_{jk}] \\
		\end{align*}
		and
		\begin{align*}
			\mathbb E [c_{ik}x_{kk}c_{jk}] = \big[(1-\delta_{ij})\beta q_k^2 p_{ik}p_{jk} + \delta_{ij} q_k^2r_k^2 p_{ik}  \big]x_{kk}
		\end{align*}
		implying 
		\begin{align*}
			\mathbb E [CXC^\top] = \alpha PQX^0QP^\top + \beta PQ \Psi^-\Psi (X) QP^\top - \beta \Psi^-(P\odot P\ \Psi(Q^2X)) + \Psi^-(P\ \Psi((QR)^2X)).
		\end{align*}
	\end{enumerate}
\end{proof}
\begin{proof}[Proof of Lemma \ref{poly}] According to the assumptions of this statement and using the previous lemma let us note the following
	\begin{enumerate}
		\item $P=P^\top$,
		\item  $\Psi X = \omega(X) \mathbf 1$ and so $\Psi^- \Psi X = \omega(X) I$, where $\omega(X)$ denotes the common diagonal element,
		\item  $X^0 = X- \Psi^-\Psi (X) = X - \omega(X) I$, 
		\item if $\beta \neq 0$ then $P\odot  P = c P$,
		\item $Q=qI$, $R=rI$.
	\end{enumerate}
	As a consequence $\mathbb E[DXC^\top], \mathbb E[DXD]$ and $\mathbb E[CXC^\top]$ consist of sums of products of the matrices $X,P,I$.
\end{proof}
\begin{proof}[Proof of Proposition \ref{prop:recursion-eigenvalues}]
	The backbone of the proof is provided by the previous lemma. The assumptions of the proposition along with Lemma \ref{poly} and Corollary \ref{PX-commute} ensure that $X_t$ and $P$ have shared eigenvectors and due to the recursion $X_t = \Phi^*(X_{t-1})$, the eigenvalues of $X_t$ are linear functions of those of $X_{t-1}$. 
\end{proof}
\begin{remark}
	Proposition \ref{prop:recursion-eigenvalues} stays true in case $P\circ P = h(P)$, where $h$ is an analytic function, in which case $c\lambda_j$ is to be replaced by $h(\lambda_j)$. 
      \end{remark}

      The next proposition, of which an extended version is to be found in \cite[Proposition 1]{gkpushsumsynchr2024}, collects the most relevant properties of $\Phi^*$.
      \begin{prop} \label{Phi*-props}
      	$\Phi^*$ exhibits the following properties:
      	\begin{enumerate}
      		\item $ \Phi^*(X)_{kl} \geq 0$ if $x_{kl} \geq 0$,
      		\item $\Phi^*(X^\top) = \Phi^*(X)^\top $,
      		\item if $X$ is PSD, then so is $\Phi^*(X)$,
      		\item if $JX = 0$, then $J\Phi^*(X)=0$ and if $XJ=0$, then $\Phi^*(X)J=0$.
      	\end{enumerate} 
      \end{prop}

      This can be used to confirm the bounding expression in \eqref{eq:errorL2bound} to be well-defined:
      \begin{corollary}
        \label{cor:PhiX0toX0}
      	Let $\mathcal X_0 = \{ X \in \mathbb R^{N\times N}: X=X^\top, XJ=0\}$. Then $\Phi^* : \mathcal X_0 \to \mathcal X_0$, whence the recursion $X_t = \Phi^*(X_{t-1})$ with $X_0 = I-J$ satisfies $X_t \in \mathcal X_0$ for all $t\geq 0$. Furthermore, since $X_0\geq 0$, we have $X_t \geq 0$ as well due to (3.) of Proposition \ref{Phi*-props}.
      \end{corollary}
      \begin{proof}[Proof of Proposition \ref{Phi*-props}]       	
      	\begin{enumerate}
      		\item Since $A$ and $X$ both consist of non-negative elements, the stated attribute is a trivial consequence of the definition of $\Phi^*$. 
      		\item Follows from:
      		\[ \Phi^*(X^\top) = \mathbb E [AX^\top A^\top] = \mathbb E[AXA^\top] ^\top = \Phi^*(X)^\top \]
      		\item Let $X$ be a PSD matrix and $w$ an arbitrary vector. Then
      		\[ w^\top\Phi^*(X)w = w^\top \mathbb E[ AXA^\top] w = \mathbb E [w^\top A X A^\top w] \geq 0 \]  
      		due to the fact that $X$ is PSD.
      		\item Let $X$ be a matrix such that $JX=0$. Then due to the linearity of the expectation we have
      		\[J\Phi^*(X) = J\mathbb E[AXA^\top] = \mathbb E[JAXA^\top] = \mathbb E[JXA^\top] = 0.\]
      		The other part is analogous and will be omitted here.
      	\end{enumerate}
      \end{proof}
      
      \section{Tools for convexity}
      \label{app:optimal}

 \begin{lemma}\label{bj>0}
 	With the notations and assumptions of Theorem \ref{generalprotocol}, we have $\forall j \  b_j\geq 0$. Furthermore if $\lambda_k > -1$ for all $k$, then $b_j > 0$ for all $j$. 
 \end{lemma}
\begin{proof}
	First let us note that the sign of $b_j$ is the same as of $(\beta- \alpha)(1-\lambda_j) - \beta c + 2r^2 =:\eta$. \\
	Case I: $\lambda_j\geq 0$. Then $(1-\lambda_j) \leq 1 $ and $\beta-\alpha \geq -1$ , hence 
	$$ \eta \geq -1 -\beta c + 2r^2. $$
	Note that due to the CSB-inequality we have $\beta c \leq r^2$, whence
	$\eta \geq 0$. \\
	Case II: $\lambda_j <0$. This means $(1-\lambda_j) \leq 2$, implying 
	$$ \eta = \beta(1-\lambda_j -c ) -\alpha(1-\lambda_j) + 2r^2 \geq -2 + 2r^2 \geq 0 $$
	since $1-\lambda_j-c \geq 0$.
	In case $\lambda_j >-1 $, then the previous inequality is strict, making $\eta>0$.
\end{proof}
\begin{prop} \label{LargestRootIneq}
  With the notations and assumptions of Theorem \ref{generalprotocol},
  let $j_0 = \arg \max_j \Delta_{j}$ and assume that $\lambda_j>-1$ for all $j$. Then
  $$\xi_1 >\Delta_{j_0}.$$
\end{prop}
\begin{proof}
	Fix $q\in(0,1)$ and let $g: \xi \mapsto f(\xi, q)$. Trivially $g$ is a polynomial of degree $N-1$ and
	\begin{align*}
		g(\Delta_{j_0}) &= -\frac{q^2}N b_{j_0} \prod_{1<j\neq j_0} (\Delta_{j_0} -\Delta_j).  
	\end{align*}
	If $j_0$ is not unique then $g(\Delta_{j_0}) = 0$, thus $\xi_1\geq \Delta_{j_0}$. In case $j_0$ is unique then $b_{j_0}\geq 0$ implies $g(\Delta_{j_0}) \leq 0$ meaning, that there must exist a root $\xi_0 \geq \Delta_{j_0}$ ($g(\xi_0)=0$) since $\lim_{\xi\to\infty} g(\xi)>0$. Trivially $\xi_1 \geq \xi_0 \geq \Delta_{j_0}$. \\
	Now let us assume that $b_{j_0} > 0$ and let $k$ denote the multiplicity of $j_0$, i.e. $k =|J| = |\{j': \max_j \Delta_j = \Delta_{j'}\}| $. Then trivially $\big.\partial_{\xi}^l f(\xi, q) \big|_{\xi=\Delta_{j_0}} = 0$ for $l< k$. For $l=k$ we have
	\begin{equation} \label{derivneg} \partial_{\xi}^{k} f(\Delta_{j_0}, q) =- \sum_{j'\in J} \frac{q^2}{N}b_j \prod_{1<j\notin J} (\Delta_{j_0} - \Delta_{j'}) < 0 
	\end{equation} 
	due to the fact that $b_{j} > 0$ (see Lemma \ref{bj>0}) and $\Delta_{j'} < \Delta_{j_0}$ for $j'\notin J$. Observe that $f$ is a polynomial in $\xi$ of the form $f=\xi^N + O(\xi^{N-1})$, hence $\lim_{\xi \to \infty} \partial_\xi^k f(\xi, q) = \infty$ for any $k<N$. Because of (\ref{derivneg}) we know that $f(\Delta_{j_0} + \varepsilon, q) <0$ for some small $\varepsilon > 0$ and since $\lim_{\xi\to \infty} f(\xi, q) = \infty$, there must exists a $\xi_1> \Delta_{j_0}$ such that $f(\xi_1, q) = 0.$

\end{proof}
\begin{remark}
	With the notations and assumptions of Theorem \ref{generalprotocol}, let $\mu := \max_j \lambda_j = \lambda_{j_0}$ and assume that $\alpha\leq 2\sqrt u /(3+(N-1)^{-1})$. Then $\max_j \Delta_j(q) = \Delta_{j_0}(q)$.
\end{remark}
\begin{proof}
	
	The proof is a series of straightforward calculations.  We need to show that the following inequality holds for any index $j$:
	\begin{align*}
		1-2q\sqrt u (1-\mu) + \alpha q^2(1-\mu)^2 \geq 1 - 2q(1-\lambda_j) + \alpha q^2(1-\lambda_j)^2.
	\end{align*}
	After rearranging terms we obtain
	\begin{align*}
		2q\sqrt u(\mu - \lambda_j) + \alpha q^2(\mu^2 - \lambda_j^2 - 2(\mu-\lambda_j)) &\geq0 \\
		q(\mu-\lambda_j)\big(2\sqrt u + \alpha q(\mu+\lambda_j -2)\big)&\geq 0.
	\end{align*}
	In case $\mu = \lambda_j$ it trivially holds, while if $\mu \geq -\frac{1}{N-1}\geq  \lambda_j\geq -1$ we have 
	$$ 2\sqrt u + \alpha q(\mu + \lambda_j -2) \geq 2\sqrt u + \alpha q(-3-\frac{1}{N-1}) \geq 0 $$
	as $q \leq 1$.
\end{proof}
\begin{lemma} 
	Let $h(\xi, q) = q^2 / (\xi- \Delta) $, where  $\Delta = \Delta(q)$ is an arbitrary smooth function. Then the determinant of the Hessian of $h$ equals
	 \begin{equation} \label{DetOfHessian}
	 	2q^4 \Delta''\cdot(\xi-\Delta)^{-5}.
	 \end{equation} 
\end{lemma}
\begin{proof}
	Since the Hessian matrix consists of the second order derivatives of $h$, we will compute these in what follows.
	\begin{align}
		\partial_\xi^2 h & = \frac{2q^2}{(\xi-\Delta)^3}, \label{eq:ddgxi}\\
		\partial_\xi\partial_q h & =  \partial_q \bigg(\frac{-q^2}{(\xi-\Delta)^2}\bigg) = \frac{-2q(\xi-\Delta)^2 - 2q^2(\xi-\Delta)  \Delta'}{(\xi - \Delta )^4} = - 2\frac{q(\xi-\Delta) + q^2 \Delta'}{(\xi-\Delta)^3}, \nonumber \\
		\partial_q^2 h & = \partial_q \bigg(\frac{2q(\xi-\Delta) + q^2\Delta'}{(\xi- \Delta)^2}\bigg) \nonumber \\ &= \frac{(2(\xi-\Delta) - 2q\Delta' + 2q\Delta ' + q^2 \Delta'')(\xi-\Delta)^2 + (2q(\xi-\Delta) + q^2 \Delta')2(\xi-\Delta)\Delta'}{(\xi-\Delta)^4} \nonumber \\
		&= \frac{2(\xi-\Delta)^2 + q^2\Delta''(\xi-\Delta) +4q(\xi-\Delta)\Delta' + 2q^2(\Delta')^2}{(\xi- \Delta)^3} \nonumber 
	\end{align}
and so
\begin{align*} \det h'' = \det \begin{bmatrix}
	\partial_\xi^2 h & \partial_q\partial_\xi h \\
	\partial_\xi \partial_q h & \partial_q^2 h
\end{bmatrix} &= \partial_\xi^2 h \partial_q^2 h - (\partial_q\partial_\xi h)^2 \\
&=\big(4q^2(\xi-\Delta)^2 + 2q^4(\xi-\Delta)\Delta'' + 8q^3(\xi-\Delta)\Delta' + 4q^4(\Delta')^2\big)/(\xi-\Delta)^6 \\
 & \quad - 4\big(q^2(\xi-\Delta)^2 + 2q^3 (\xi-\Delta)\Delta' + q^4 (\Delta')^2\big)/(\xi-\Delta)^6 \\
 &= \frac{2q^4(\xi-\Delta) \Delta''}{(\xi-\Delta)^6} = 2q^4 \Delta'' (\xi-\Delta)^{-5}.
\end{align*}
\end{proof}

We now apply the generic lemma above for our more specific context.

\begin{corollary} \label{PosDefHessian}
	Let $\Delta(q) = 1- 2q(1-\lambda) + \alpha q^2(1-\lambda)^2$ in line with the function used in Theorem \ref{generalprotocol}. Then the previous determinant (\ref{DetOfHessian}) takes the form
	$$ 4q^4\alpha(1-\lambda)^2 (\xi-\Delta)^{-5}.$$
      \end{corollary}

\begin{corollary}
	The Hessian of $h$ defined above is positive semi-definite at $(\xi, q)$ if and only if $\alpha\geq 0$ and $\xi\geq \Delta(q)$.
      \end{corollary}
      \begin{proof}
        The matrix is positive semi-definite if all the determinants of the principal minors are non-negative. Besides the overall determinant considered in the corollary above, the first entry is also non-negative as is expressed in \ref{eq:ddgxi}, given the current inequalities imposed.
      \end{proof}

  \begin{lemma}
  	Let $g(\xi,q)$ be a convex, twice continuously differentiable function with non-vanishing derivative $\partial_\xi g$, and let $\xi:[0,1] \to [0,1]$ be the solution to $g(\xi,q)=0$. If $\partial_\xi g \leq 0$ ($\partial_\xi g \geq 0$) in some open neighborhood of $(\xi(q_0), q_0)$, then $\xi(q)$ is convex (concave resp.) in an open neighborhood of $q_0$. 
  \end{lemma}
\begin{proof}
	Calculating the second derivative of the equation $g(\xi(q), q) = 0$ yields 
	\begin{align*}
		0 = \partial_\xi^2 g\cdot (\xi')^2 + 2 \partial^2_{q\xi} g\cdot \xi' + \partial_q^2 g + \partial_\xi g\cdot \xi''
	\end{align*}
which implies
$$ \xi''(q_0) = -\frac{1}{\partial_{\xi}g(\xi(q_0),q_0)} \begin{bmatrix}
	\xi'(q_0)& 1 
\end{bmatrix} H_g(\xi(q_0), q_0) \begin{bmatrix}
\xi'(q_0) \\ 1
\end{bmatrix} $$
where $H_g$ denotes the Hessian matrix of $g$. Since $H_g(\xi(q_0), q_0)\geq 0$ the sign of $\xi''(q_0)$ depends on that of $\partial_\xi g(\xi(q_0),q_0)$.
\end{proof}
      \begin{proof}[Proof of Proposition \ref{prop:lines_xidiff}]
First we will calculate $\xi_1'(q)$ using the implicit function theorem. Since $$ f(\xi_1(q), q) =0 \quad \forall q\in(0,1) $$
we have
\begin{equation*}
	0 = (\partial_q f(\xi_1(q),q)) = \partial_\xi f (\xi_1(q), q) \xi_1'(q) + \partial_q f(\xi_1(q), q)
\end{equation*}
from which we conclude
\[ \xi_1'(q) = -\frac{\partial_q f(\xi_1(q),q)}{\partial_\xi f(\xi_1(q),q)} .\]
Now let $\pi(\xi, q) = \prod_{j>1}(\xi - \Delta_j(q))$ and $\sigma(\xi,q) = 1- q^2/N\sum_{j>1}b_j/(\xi-\Delta_j(q))$, with which we can write $f = \pi \sigma$. Easy calculation shows 
\begin{align*}
	\partial_\xi f(\xi, q) =\partial_\xi \pi(\xi,q) \bigg(1- \frac{q^2}{N}\sum_{j>1} \frac{b_j}{\xi-\Delta_j(q)}\bigg) + \pi(\xi,q) \frac{q^2}{N}\sum_{j>1} \frac{b_j}{(\xi-\Delta_j(q))^2}  
\end{align*}
and 
\begin{align*}
	\partial_q f(\xi, q) &= \partial_q \pi(\xi,q)\bigg(1- \frac{q^2}{N}\sum_{j>1}\frac{b_j}{\xi-\Delta_j(q)}\bigg) + \pi(\xi,q)\bigg(-\frac{2q}{N}\sum_{j>1}\frac{b_j}{\xi - \Delta_j(q)} - \frac{q^2}{N}\sum_{j>1} \frac{b_j \partial_q\Delta_j(q)}{(\xi-\Delta_j(q))^2}\bigg).
\end{align*}
We know from Proposition \ref{LargestRootIneq} that $\xi_1(q) > \max_j \Delta_j$, so $\pi(\xi_1, q) \neq 0$, thus $\sigma(\xi_1,q) = 0$ must hold, alas
 \begin{equation} \label{XiDiff}\xi_1'(q) =- \frac{\partial_q\sigma(\xi_1, q)}{\partial_\xi \sigma(\xi_1,q)} = \frac{2q\sum_{j>1}b_j/(\xi_1-\Delta_j(q)) + q^2\sum_{j>1}  b_j\partial_q\Delta_j(q) / (\xi_1-\Delta_j(q))^2 ) }{ q^2 \sum_{j>1} b_j / (\xi_1-\Delta_j(q))^2} 
 \end{equation}
where $\Delta_j'(q) = -2(1-\lambda_j) + 2\alpha q (1-\lambda_j)^2. $ \\
Now let us calculate $\xi_1'(0)$, in which case trivially $\xi_1(0) = \max_j \Delta_j(0) = 1$. Denote by $J=\{j': \max_j\Delta_j(0)= \Delta_{j'}\}$. Note that $\partial_q f(1,0) = 0 = \partial_q f(1,0)$, hence we cannot use the implicit function theorem directly. 
First we observe that due to $\xi_1 : [0,1] \to \mathbb R$ being continuous and convex on $(0,1)$, $\xi_1$ is convex on $[0,1]$ as well, furthermore its right-derivatives are continuous from the right, i.e.
$$ \xi'_1(0) = \lim_{q\to 0+} \frac{\xi_1(q) - \xi_1(0)}{q} = \lim_{q\to 0+} \xi'_{1,+}(q). $$
We have shown before that $\xi_1'$ exist for every $q\in(0,1)$, so 
\begin{align*} \xi_1'(0) &= \lim_{q\to 0+} \frac{\partial_q \sigma(\xi_1(q),q)}{\partial_\xi \sigma(\xi_1(q), q)} = \frac{\sum_{j\in J} b_j \Delta_{j}'(0)}{\sum_{j\in J} b_j} = \Delta_{j_0}'(0).
\end{align*}
Trivially $\Delta_{j_0}'(0) = -2(1-\lambda_j)$. 
Due to convexity, $\xi_1(q) \geq \xi_1'(0) q +1$ and $\xi_1(q) \geq \xi_1'(1)(q-1) + 1 $ has to hold for any $q\in(0,1)$.	
\end{proof}

\end{document}